\documentclass[12pt,a4paper]{article}

\usepackage{amssymb}

\usepackage{amsmath}

\newtheorem{theorem}{Theorem}[section]
\newtheorem{definition}[theorem]{Definition}
\newtheorem{lemma}[theorem]{Lemma}

\begin{document}
\title{Gr\"{o}bner-Shirshov basis for HNN extensions of groups and for the alternating group \footnote
{The research is supported by the National Natural Science
Foundation of China (Grant No.10771077) and the Natural Science
Foundation of Guangdong Province (Grant No.06025062).}}
\author{
Yuqun Chen and Chanyan Zhong  \\
\\
{\small \ School of Mathematical Sciences}\\
{\small \ South China Normal University}\\
{\small \ Guangzhou 510631}\\
{\small \ P. R. China}\\
{\small \ yqchen@scnu.edu.cn} \\
{\small \ zhongchanyan@tom.com}}
\date{}
\maketitle \noindent\textbf{Abstract:} In this paper, we generalize
the Shirshov's Composition Lemma by replacing the monomial order for
others. By using Gr\"{o}bner-Shirshov bases, the normal forms of HNN
extension of a group and the alternating group are obtained.

\noindent \textbf{Key words: }Gr\"{o}bner-Shirshov basis,
 normal form, HNN extension, alternating group

\noindent {\bf AMS} Mathematics Subject Classification(2000): 20E06,
16S15, 13P10

\section{Preliminaries}

It is known that in the Gr\"{o}bner-Shirshov basis theory, the
Shirshov's Composition-Diamond Lemma \cite{s} plays an important
role. In the Composition-Diamond Lemma, the order is asked to be
monomial. In this paper, we generalize the Composition-Diamond Lemma
by replacing the monomial order for others. From this result, we
show by direct calculations of compositions, that the presentation
of the HNN extension is a Gr\"{o}bner-Shirshov basis under an
appropriate ordering of group words, in which the order is not
monomial order. By the generalized composition lemma, we immediately
obtain the Normal Form Theorem for HNN extensions. In fact, this is
the first time to find a Gr\"{o}bner-Shirshov basis by using a
non-monomial order.

HNN extensions of groups were first invented by
Higman-Neumann-Neumann in 1949 (\cite{hnn}) and independently by P.
S. Novikov in 1952 (see \cite{n52}, \cite{n54}, \cite{n55}). The
Normal Form Theorem for certain HNN extensions of groups was first
established by L. A. Bokut (see \cite{b66}, \cite{b67}, \cite{b68}
and see also K. Kalorkoti \cite{k}). In general, the Normal Form
Theorem for HNN extensions of groups was proved in the text of
Lyndon-Schupp \cite{Ly}.

For the alternating group $A_n$, a presentation was given in the
monograph of N. Jacobson (see \cite{bu}, p71). However, we still do
not know what is the normal form for $A_n$. In this paper, we give
the normal form theorem of $A_n$ with respect to the above
presentation.
\\

We first cite some concepts and results from the literature. Let $k$
be a field, $k \langle X\rangle$ the free associative algebra over
$k$ generated by $X$ and $ X^{*}$ the free monoid generated by $X$,
where the empty word is the identity which is denoted by 1. For a
word $w\in X^*$, we denote the length of $w$ by $deg(w)$. Let $X^*$
be a well ordered set. Let $f=\sum_{a\in X^*}f(a)a\in k\langle
X\rangle$, where $f(a)\in k$, with the leading word $\bar{f}$. We
say that $f$ is monic if $\bar{f}$ has coefficient 1. We denote
$suppf=\{a\in X^*|f(a)\neq 0\}$.


\begin{definition} (\cite{s}, see also \cite{b72}, \cite{b76}) \
Let $f$ and $g$ be two monic polynomials in $k\langle X\rangle$.
Then, there are two kinds of compositions:

$(1)$ If \ $w$ is a word such that $w=\bar{f}b=a\bar{g}$ for some
$a,b\in X^*$ with deg$(\bar{f})+$deg$(\bar{g})>$deg$(w)$, then the
polynomial
 $(f,g)_w=fb-ag$ is called the intersection composition of $f$ and
$g$ with respect to $w$.

$(2)$ If  $w=\bar{f}=a\bar{g}b$ for some $a,b\in X^*$, then the
polynomial $(f,g)_w=f - agb$ is called the inclusion composition of
$f$ and $g$ with respect to $w$.

In the above case, the transformation $f\mapsto(f,g)_w=f-agb$ is
called the elmination of the leading word (ELW) of $g$ in $f$.
\end{definition}

\begin{definition}(\cite{b72}, \cite{b76}, cf. \cite{s})
Let $S\subseteq k\langle X\rangle$ and $<$ a well order on $X^*$.
Then the composition $(f,g)_w$ is called trivial modulo $S$ if
$(f,g)_w=\sum\alpha_i a_i s_i b_i$, where each $\alpha_i\in k$,
$a_i,b_i\in X^{*}$ and $\overline{a_i s_i b_i}<w$. If this is the
case, then we write
$$
(f,g)_w\equiv0\quad mod(S,w)
$$
In general, for $p,q\in k\langle X\rangle$, we write
$$
p\equiv q\quad mod (S,w)
$$
which means that $p-q=\sum\alpha_i a_i s_i b_i $, where each
$\alpha_i\in k,a_i,b_i\in X^{*}$ and $\overline{a_i s_i b_i}<w$.
\end{definition}

\begin{definition} (\cite{b72}, \cite{b76}, cf. \cite{s}) \
We call the set $S$ with respect to the well order $``<"$a
Gr\"{o}bner-Shirshov set (basis) in $k\langle X\rangle$ if any
composition of polynomials in $S$ is trivial modulo $S$.
\end{definition}

\noindent{\bf Remark}: Usually, in the definition of the
Gr\"{o}bner-Shirshov basis, the order is asked to be monomial.

A well order $``<"$ on $X^*$ is monomial if it is compatible with
the multiplication of words, that is, for $u, v\in X^*$, we have
$$
u > v \Rightarrow w_{1}uw_{2} > w_{1}vw_{2},  \ for \  all \
 w_{1}, \ w_{2}\in  X^*.
$$

The following lemma was proved by Shirshov \cite{s} for the free Lie
algebras (with deg-lex ordering) in 1962 (see also Bokut
\cite{b72}). In 1976, Bokut \cite{b76} specialized the approach of
Shirshov to associative algebras (see also Bergman \cite{b}). For
commutative polynomials, this lemma is known as the Buchberger's
Theorem (see \cite{bu65}), published in \cite{bu70}.

\begin{lemma}\label{l1}
(Composition-Diamond Lemma) \ Let $A=k\langle X|S\rangle$ and $``<"$
a monomial order on $X^*$. Then the following statements are
equivalent:
\begin{enumerate}
\item[(i)] \ $S$ is a Gr\"{o}bner-Shirshov basis.
\item[(ii)] \ For any $f\in k\langle X\rangle, \
0\neq f\in Ideal(S)\Rightarrow \bar{f}=a\bar{s}b$ for some $s\in S,
\ a,b\in X^*$.
\item[(iii)] \
The set
$$
Red(S)=\{u\in X^*|u\neq a\bar{s}b,s\in S,a,b\in X^*\}
$$
is a linear basis of the algebra $A$.
\end{enumerate}
\end{lemma}

If a subset $S$ of $k\langle X\rangle$ is not a Gr\"{o}bner-Shirshov
basis, then we can add to $S$ all nontrivial compositions of
polynomials of $S$, and by continuing this process (maybe
infinitely) many times, we eventually obtain a Gr\"{o}bner-Shirshov
basis $S^{comp}$. Such a process is called the Shirshov algorithm.

If $S$ is a set of ``semigroup relations" (that is, the polynomials
of the form $u-v$, where $u,v\in X^{*}$), then any nontrivial
composition will have the same form. As a result, the set $S^{comp}$
also consists of semigroup relations.

Let $A=sgp\langle X|S\rangle$ be a semigroup presentation. Then $S$
is a subset of $k\langle X\rangle$  and hence one can find a
Gr\"{o}bner-Shirshov basis $S^{comp}$. The last set does not depend
on $k$, and as mentioned before, it consists of semigroup relations.
We will call $S^{comp}$ a Gr\"{o}bner-shirshov basis of $A$. This is
the same as a Gr\"{o}bner-shirshov basis of the semigroup algebra
$kA=k\langle X|S\rangle$.

\section{Generalized Composition-Diamond Lemma}

In this section, we generalize the Composition-Diamond Lemma which
is useful in the sequel, by replacing the monomial order for others.
The proof of the following lemma is essentially the same as in
\cite{bf}. For the sake of convenience, we give the details.

\begin{lemma}\label{l2}
(Generalized Composition-Diamond Lemma) \ Let $S\subseteq k\langle
X\rangle, \ A=k\langle X|S\rangle$ and $``<"$ a well order on $X^*$
such that
\begin{enumerate}
\item[(A)] \ $\overline{asb}=a\bar{s}b$ for any $a,b\in X^*, \ s\in S$;
\item[(B)] \ for each composition $(s_1,s_2)_w$ in $S$, there exists a presentation
$$
(s_1,s_2)_w=\sum_{i}\alpha_{i}a_it_ib_i, \ a_i\bar{t_i}b_i<w, \ \
\mbox{ where } \ t_i\in S, \ a_i,b_i\in X^*, \ \alpha_{i}\in k
$$
 such that for any $c,d\in X^*$, we have
\begin{eqnarray}
ca_i\bar{t_i}b_id<cwd
\end{eqnarray}
\end{enumerate}
Then, the following statements hold.
\begin{enumerate}
\item[(i)] \ $S$ is a Gr\"{o}bner-Shirshov basis.
\item[(ii)] \ For any $f\in k\langle X\rangle, \
0\neq f\in Ideal(S)\Rightarrow \bar{f}=a\bar{s}b$ for some $s\in S,
\ a,b\in X^*$.
\item[(iii)] \
The set
$$
Red(S)=\{u\in X^*|u\neq a\bar{s}b,s\in S,a,b\in X^*\}
$$
is a linear basis of the algebra $A$.
\end{enumerate}
\end{lemma}

\noindent{\bf Proof} \ (i) is clear. Now, we prove (ii). Let
$$ f=\sum_{i=1}^n\alpha_ia_is_ib_i,
\ \alpha_i\in k, \ s_i\in S, \ a_i,b_i\in X^*$$ Assume that
$$w_i=a_i\bar{s_i}b_i, \ w_1=w_2=\cdots=w_l>w_{l+1}\cdots$$ We will
use the induction on $l$ and $w_1$ to prove that
$\bar{f}=a\bar{s}b$, for some $s\in S \ \mbox{and} \ a,b\in X^*$.

If $l=1$, then by (A), $\bar{f}=a_1\bar{s_1}b_1$ and hence the
result holds. Assume that $l\geq 2$. Then
$$
f=(\alpha_1+\alpha_2)a_1s_1b_1+\alpha_2(a_2s_2b_2-a_1s_1b_1)+\cdots
$$
For $w_1=w_2$, there are three cases to consider.

Case 1. Assume that $b_1=b\bar{s_2}b_2$ and $a_2=a_1\bar{s_1}b$.
Then we have
$$
a_2s_2b_2-a_1s_1b_1=a_1s_1b(s_2-\bar{s_2})b_2-a_1(s_1-\bar{s_1})bs_2b_2.
$$
For any $t\in supp(s_2-\bar{s_2})$, by (A),
$\overline{a_1s_1btb_2}=a_1\bar{s_1}btb_2<a_1\bar{s_1}b\bar{s_2}b_2=w_1$
and similarly, we have $\overline{a_1t_1bs_2b_2}<w_1$, for any
$t_1\in supp(s_1-\bar{s_1})$.

Case 2.  Assume that $b_1=bb_2, \ a_2=a_1a, \ \bar{s_1}b=a\bar{s_2}$
and $deg\bar{s_1}+deg\bar{s_2}>deg(a\bar{s_1})$. Then
$$
a_2s_2b_2-a_1s_1b_1=a_1(as_2-s_1b)b_2.
$$
By (B), there exist $\beta_j\in k, \ u_j,v_j\in X^*, \ t_j\in S$
such that $u_j\bar{t_j}v_j<w=\bar{s_1}b, \
s_1b-as_2=\sum_j\beta_ju_jt_jv_j$ and
$a_1u_j\bar{t_j}v_jb_2<a_1\bar{s_1}bb_2$. Now, by (A), for any $j$,
we have
$\overline{a_1u_jt_jv_jb_2}=a_1u_j\bar{t_j}v_jb_2<a_1\bar{s_1}bb_2=w_1$.

Case 3.  Assume that $b_2=bb_1, \ a_2=a_1a$ and
$\bar{s_1}=a\bar{s_2}b$. Then
$$
a_2s_2b_2-a_1s_1b_1=a_1(as_2b-s_1)b_1.
$$
by (A) and (B), there exist $\beta_j\in k, \ u_j,v_j\in X^*, \
t_j\in S$ such that $u_j\bar{t_j}v_j<w=\bar{s_1}, \
s_1-as_2b=\sum_j\beta_ju_jt_jv_j$ and for any $j, \ \mbox{we have} \
\overline{a_1u_jt_jv_jb_1}=a_1u_j\bar{t_j}v_jb_1<a_1\bar{s_1}b_1=w_1$.

(iii) follows from (ii).  \ \ $\square$

\section{Normal form for HNN extension of group}

In this section, by using the Gr\"{o}bner-Shirshov basis, we provide
a new proof of the normal form theorem of HNN extension of a group.

\begin{definition} (\cite{hnn},\cite{n52},\cite{n54},\cite{n55}) \
Let $G$ be a group and let $A$ and $B$ be
subgroups of $G$ with $\phi:A\rightarrow B$ an isomorphism. Then the
HNN extension of $G$ relative to $A\ ,B$ and $\phi$ is the group
$$
{\cal G}=gp\langle G,t;t^{-1}at=\phi(a),a\in A\rangle.
$$
\end{definition}
Let $G=G_1\cup\{1\}$, where $G_1=G\setminus\{1\}=
\{g_\alpha;\alpha\in\Lambda\}$, and let
$$
G/A=\{g_i A;i\in I\}, \ G/B=\{h_j B;j\in J\},
$$
where $\{g_i;i\in I\}$ and $\{h_j;j\in J\}$ are the coset
representatives of $A$ and $B$ in $G$, respectively. We assume that
all sets $\Lambda, \ I, \ J$ are well ordered and so are the sets
$\{g_\alpha;\alpha\in\Lambda\}, \ \{g_i;i\in I\}, \ \{h_j;j\in J\}$.
Then we get a new presentation of the group ${\cal G}$ as a
semigroup:
\begin{eqnarray*}
{\cal G}=sgp\langle G_1,t,t^{-1}&;&g g^{'}=[g g^{'}],g t=g_{_A}
t\phi(a_g),gt^{-1}=g_{_B}
t^{-1}\phi^{-1}(b_g),\\
&&t^{\varepsilon}t^{-\varepsilon}=1, \ g,g'\in G_1, \
\varepsilon=\pm1\rangle,
\end{eqnarray*}
where $[g g^{'}]\in G; \ g=g_{_A} a_g, \ g_{_A}\in\{g_i;i\in I\}, \
g\neq g_{_A}, \ a_g\in A; \ g=g_{_B} b_g, \ g_{_B}\in\{h_j;j\in J\},
\ g\neq g_{_B}, \ b_g\in B$.

Now, we order the set $G$ in three different ways:
\begin{enumerate}
\item[(1)]\ Let $1<g_{_\alpha}<g_{_\beta}<\cdots \ (\alpha<\beta)$ be a well order of $G$.
Then we denote this order by $(G,>)$ and call it an absolute order.
\item[(2)]\ For any $g,g^{'}\in G$, suppose that $g=g_{_A}
a_g,g^{'}=g^{'}_{_A} a_{g^{'}}$. Then $g>_{_A}g^{'}$ if and only if
$(g_{_A},a_g)>(g^{'}_{_A},a_g^{'})$ is ordered lexicographically
(elements $g_{_A},g^{'}_{_A}$ by $I$, elements $a_g,a_{g^{'}}$ by
(1)). We denote this order by $(G,>_{_A})$ and call it the
$A$-order. In particular, if $g\neq g_{_A}$, then $g>_{_A}g_{_A}, \
\mbox{for} \ (g_{_A},a_g)>(g_{_A},1), \ a_g\neq1$.
 \item[(3)]\ For any $g,g^{'}\in G$, suppose that
$g=g_{_B} b_g,g^{'}=g'_{_B} b_{g^{'}}$. Then $g>_{_B}g^{'}$ if and
only if $(g_{_B},b_g)>(g'_{_B},b_g^{'})$ is ordered
lexicographically (elements $g_{_B},g^{'}_{_B}$ by $J$, elements
$b_g,b_{g^{'}}$ by (1)). We denote this order by $(G,>_{_B})$ and
call it the $B$-order.
\end{enumerate}

Then we order the set $G^{*}_1$ in three different ways too:
\begin{enumerate}
\item[(1)]\ The absolute order $({G_1}^{*},\leq)$ is deg-lex order, to compare words $g_1\cdots
g_n \ (n\geq0)$ first by length and then lexicographically using
absolute order of $G_1$.
\item[(2)]\ The $A$-order $({G_1}^{*},\leq_{_A})$ is ${deg-lex}_{_A}$ order,
firstly to compare words $g_1\cdots
g_n \ (n\geq0)$ by length, secondly for $n\geq1, \
g_1,\cdots,g_{n-1}$ lexicographically by absolute order, and finally
the last elements $g_n$ by $A$-order.
 \item[(3)]\ The $B$-order $({G_1}^{*},\leq {_B})$ is similar to
 (2) replacing $>_{_A} \ \mbox{by} \ >_{_B}$.
\end{enumerate}

Each element in $\{G_1\dot{\cup}\{t,t^{-1}\}\}^*$ has a unique form
$u=u_1 t^{{\varepsilon}_1} u_2 t^{{\varepsilon}_2}\cdots
u_kt^{{\varepsilon}_k} u_{k+1}$, where each ${u_i}\in G_1^*,
{\varepsilon_i}=\pm1, \ k\geq 0$. Suppose that $v=v_1
t^{{\delta}_1}v_2\cdots
v_lt^{{\delta}_l}v_{l+1}\in\{G_1\dot{\cup}\{t,t^{-1}\}\}^*$.

Then
$$wt(u)=(k,t^{\varepsilon_1},\cdots,t^{\varepsilon_k},u_1,\cdots,u_k,u_{k+1})$$
$$wt(v)=(l,t^{\delta_1},\cdots,t^{\delta_l},v_1,\cdots,v_l,v_{l+1})$$
We define $u\succ v$ if $wt(u)>wt(v)$ lexicographically, using the
order of natural numbers and the following orders:
\begin{enumerate}
\item[(a)]\ $t>t^{-1}$
\item[(b)]\ $u_i>_{_A}v_i \ \mbox{if} \ \varepsilon_i=1, \ 1\leq i\leq k$
\item[(c)]\ $u_i>_{_B}v_i \ \mbox{if} \ \varepsilon_i=-1, \ 1\leq i\leq k$
\item[(d)]\ $u_{k+1}> v_{l+1} \ (k=l)$, the absolute order of ${G_1}^{*}$
 \end{enumerate}

Now, we can easily verify the following lemma.

\begin{lemma}\label{l3}
Let the order $\succ$ on $\{G_1\dot{\cup}\{t,t^{-1}\}\}^*$ be
defined as above. Then the order $\succ$ is a well order but not
monomial, for example, $g\succ g^{'}$ does not necessarily imply
that $gt\succ g^{'}t$.
\end{lemma}

Equipping with the above notation, we have the following lemma.

\begin{lemma}\label{l4}
Let $X=G_1\dot{\cup}\{t,t^{-1}\}$. Suppose that the order $\succ$ on
$X^*$ is defined as above and $S=\{gg^{'}-[gg^{'}], \
gt-g_{_A}t\phi(a_{g}), \ gt^{-1}-g_{_B}t^{-1}\phi^{-1}(b_{g}), \
t^{\varepsilon}t^{-\varepsilon}-1 \ | \ g,g'\in G_1, \
\varepsilon=\pm1\}$ is as above too. Then $S$  satisfies conditions
(A)-(B) in Lemma \ref{l2}.
\end{lemma}
\textbf{Proof} \ For any $c,d\in X^*$, suppose that $c=c_1
t^{\varepsilon_1}\cdots c_{n}t^{\varepsilon_n} c_{n+1}, \ d=d_1
t^{\delta_1}\cdots d_{m}t^{\delta_m} d_{m+1}, \ c_i,d_j\in G_1^*, \
\varepsilon_i,\delta_j=\pm1$. We firstly check (A). Then there are
four cases to consider. For example, the second case is for
polynomials $gt-g_{_A}t\phi(a_{g}), \ g\neq g_{_A}$. We need to
prove that $cgtd\succ cg {_A}t\phi(a_{g})d$ for any $c,d$.

Since
\begin{eqnarray*}
wt(cgtd)&=&(n+m+1,t^{\varepsilon_1},\cdots,t^{\varepsilon_n},t,
t^{\delta_1},\cdots,t^{\delta_m},c_1,\cdots,
c_{n},c_{n+1}g,\\
&&d_1,\cdots,d_{m},d_{m+1}),\\
wt(cg_{_A}t\phi(a_{g})d)&=&(n+m+1,t^{\varepsilon_1},\cdots,t^{\varepsilon_n},t,
t^{\delta_1},\cdots,t^{\delta_m},c_1,\cdots,
c_{n},c_{n+1}g_{_A},\\
&&\phi(a_{g})d_1,\cdots,d_{m},d_{m+1})
\end{eqnarray*}
and $c_{n+1}g>_{_A}c_{n+1}g_{_A}$ (since $g>_{_A}g_{_A}$), we have
$cgtd\succ cg {_A}t\phi(a_{g})d$.

We secondly check that (B) holds in Lemma \ref{l2}. By noting that
there is no inclusion compositions in $S$, we need only to consider
the cases of intersection compositions. For any $a,b\in X^*, \
s_1,s_2\in S$, suppose that $a\bar{s_1}=\bar{s_2}b$ with
$deg\bar{s_1}+deg\bar{s_2}>deg(a\bar{s_1})$. Then, we consider the
following cases:
$$
w=gg^{'}g^{''}, \ gg^{'}t \ (g'\neq g'_{_A}), \ gg^{'}t^{-1} \
(g'\neq g'_{_B}), \ gt^{\varepsilon}t^{-\varepsilon}, \
t^{\varepsilon}t^{-\varepsilon}t^{\varepsilon} \ (\varepsilon=\pm1).
$$
For example, the second case is as follows:

 \ Let $\bar{s_1}=gg^{'}, \ \bar{s_2}=g't, \ w=gg^{'}t$. Then,
by noting that
$gg^{'}=[gg^{'}]_{_A}a_{[gg^{'}]}=[gg^{'}_{_A}]_{_A}a_{[gg^{'}_{_A}]}a_{g'}$
implies that $[gg^{'}]_{_A}=[gg^{'}_{_A}]_{_A}$ and
$a_{[gg^{'}]}=a_{[gg^{'}_{_A}]}a_{g'}$, we know that
$(s_2,s_1)_w=[gg^{'}]t-g{g^{'}}_{_A}t\phi(a_{g'})=([gg^{'}]t-[gg^{'}]_{_A}t\phi(a_{[gg^{'}]}))-
(gg^{'}_{_A}-[gg^{'}_{_A}])t\phi(a_{g'})-([g{g^{'}}_{_A}]t\phi(a_{g'})-
[gg^{'}_{_A}]_{_A}t\phi(a_{[gg^{'}_{_A}]}a_{g'}))$. We denote
$s'_1=[gg^{'}]t-[gg^{'}]_{_A}t\phi(a_{[gg^{'}]}), \
s'_2=gg^{'}_{_A}-[gg^{'}_{_A}], \ b_2=t\phi(a_{g'}), \
s'_3=[g{g^{'}}_{_A}]t\phi(a_{g'})-
[gg^{'}_{_A}]_{_A}t\phi(a_{[gg^{'}_{_A}]}a_{g'})$. Clearly, $s'_1, \
s'_2, \ s'_3\in S \ \mbox{or} \ \{0\}$ (if
$[gg^{'}]=[gg^{'}]_{_A}$). Since
\begin{eqnarray*}
wt(c\bar{s'_1}d)&=&(n+m+1,t^{\varepsilon_1},\cdots,t^{\varepsilon_n},t,
t^{\delta_1},\cdots,t^{\delta_m},c_1,\cdots,
c_{n},c_{n+1}[gg^{'}],\\
&&d_1,\cdots,d_{m},d_{m+1}),\\
wt(c\bar{s'_2}b_2
d)&=&(n+m+1,t^{\varepsilon_1},\cdots,t^{\varepsilon_n},t,
t^{\delta_1},\cdots,t^{\delta_m},c_1,\cdots,
c_{n},c_{n+1}gg^{'}_{_A},\\
&&\phi(a_{g'})d_1,\cdots,d_{m},d_{m+1}),\\
wt(c\bar{s'_3}d)&=&(n+m+1,t^{\varepsilon_1},\cdots,t^{\varepsilon_n},t,
t^{\delta_1},\cdots,t^{\delta_m},c_1,\cdots,
c_{n},c_{n+1}[gg^{'}_{_A}],\\
&&\phi(a_{g'})d_1,\cdots,d_{m},d_{m+1}),\\
wt(cwd)&=&(n+m+1,t^{\varepsilon_1},\cdots,t^{\varepsilon_n},t,
t^{\delta_1},\cdots,t^{\delta_m},c_1,\cdots,
c_{n},c_{n+1}gg^{'},\\
&&d_1,\cdots,d_{m},d_{m+1})
\end{eqnarray*}
and $c_{n+1}gg^{'}>_{_A}c_{n+1}[gg^{'}], \ c_{n+1}gg^{'}_{_A}, \
c_{n+1}[gg^{'}_{_A}]$, we have $cwd\succ c\bar{s'_1}d, \
c\bar{s'_2}b_2 d, \ c\bar{s'_3}d$. All other cases are treated the
same. The proof is finished. \ \ $\square$
\\

Now, by Lemma \ref{l2} and Lemma \ref{l4}, we obtain the following
theorem.
\begin{theorem}\label{t1}
A Gr\"{o}bner-Shirshov basis of HNN extention ${\cal G}=gp\langle
G,t;t^{-1}at=\phi(a),a\in A\rangle$ consists of the following
relations:

\begin{enumerate}
\item[(3.1)] \ $gg^{'}=[gg^{'}]$;
\item[(3.2)] \
$gt=g_{_A}t\phi(a_{g})$, where $g=g_{_A}a_{g}, \
g_{_A}\in\{g_{i},i\in I\}, \ a_{g}\in A$;
\item[(3.3)] \
$gt^{-1}=g_{_B}t^{-1}\phi^{-1}(b_{g})$, where $g=g_{_B}b_{g}, \
g_{_B}\in\{h_{j},j\in J\}, \ b_{g}\in B$;
\item[(3.4)] \ $tt^{-1}=1, \ t^{-1}t=1$.
\end{enumerate}
where $\{g_i;i\in I\}$ and $\{h_j;j\in J\}$ are the coset
representatives of $A$ and $B=\phi(A)$ in $G$, respectively.
\end{theorem}

\ \

Thus $S$ is a Gr\"{o}bner-Shirshov basis of the algebra $k\langle
G_1\dot{\cup}\{t,t^{-1}\}\rangle$ with
$$
Red(S)=\{u\in
\{G_1\dot{\cup}\{t,t^{-1}\}\}^*|u\neq a\bar{s}b,s\in S,{a,b}\in
\{G_1\dot{\cup}\{t,t^{-1}\}\}^*\}
$$
which is the normal form of the HNN extension ${\cal G}$.  From this
result, the normal form theorem for HNN Extension of the Group $G$
easily follows.

\begin{theorem} \label{th5}
(The Normal Form Theorem for HNN Extension, \cite{Ly}, Theorem
4.2.1) Let  ${\cal G}=gp\langle G,t;t^{-1}at=\phi(a),a\in A\rangle$
be an HNN extension of group $G$. If $\{g_{i};i\in I\}, \ \mbox{and}
\ \{h_{j};j\in J\}$ are the sets of representatives of the left
cosets of $A$ and $B=\phi(A)$ in $G$, respectively, then every
element $w$ of ${\cal G}$ has a unique representation $w=g_1
t^{\varepsilon_1} \cdots g_n t^{\varepsilon_n} g_{n+1} \ (n\geq 0, \
\varepsilon_l=\pm1)$, where, for $1\leq l\leq n$, the following
conditions are satisfied:
\begin{enumerate}
\item[(1)] \ if $\varepsilon_{l}=1$, then $g_{l}\in \{g_{i};i\in
I\}$,
\item[(2)] \ if $\varepsilon_{l}=-1$, then $g_{l}\in
\{h_{j};j\in J\}$,
\item[(3)] \ there does not exist subwords $tt^{-1} \ and \ t^{-1}t$,
\item[(4)] $g_{n+1}$ is an arbitrary element of $G$.
\end{enumerate}
\end{theorem}
\ \

\noindent{\bf Remark}: In  Theorem 4.2.1 of \cite{Ly}, the right
cosets were considered. We notice that the above Theorem \ref{th5}
is essentially the same as  Theorem 4.2.1 in \cite{Ly}.

\section{Normal form for alternating group}

In this section, we first find a Gr\"{o}ner-Shirshov basis for the
alternating group $A_n$ and then we give the normal form theorem of
$A_n$.

Let $S_{n}$ be the group of the permutations of $\{1,2,\cdots, n\}$.
Then the subset $A_{n}$ of all even permutations in $S_n$ is a
normal subgroup of $S_{n}$. We call $A_{n}$ the alternating group of
degree $n$. The following presentation of $A_{n}$ was given in the
monograph of N. Jacobson (see \cite{bu}, p71):
\begin{eqnarray*}
A_{n}=gp\langle x_{i} \ (1\leq i\leq n-2)&;&x_{1}^{3}=1, \
(x_{i-1}x_{i})^3=x_{i}^{2}=1 \ (2\leq i\leq n-2),\\
&&(x_{i}x_{j})^{2}=1 \ (1\leq i<j-1,j\leq n-2)\rangle
\end{eqnarray*}
where $x_{i}=(12)((i+1)(i+2)), \ i=1,2,\cdots, n-2, \ (ij)$ the
transposition. We now give a presentation of the group $A_{n}$ as a
semigroup:
\begin{eqnarray*}
A_{n}&=&sgp\langle x_{1}^{-1}, \ x_{i} \ (1\leq i\leq n-2); \
 x_{1}x_{1}^{-1}=
x_{1}^{-1}x_{1}=x_{1}^{3}=1, \ (x_{i-1}x_{i})^3=x_{i}^{2}=1\\
&&(2\leq i\leq n-2), \ (x_{i}x_{j})^{2}=1 \ (1\leq i<j-1,j\leq
n-2)\rangle
\end{eqnarray*}
We now order the generators in the following way:
$$x_{1}^{-1}<x_{1}<x_{2}<\cdots <x_{n-2}$$
Let $X=\{x_1^{-1},x_1,x_2,\cdots,x_{n-2}\}$. Then, with the above
notations, we can order the words of $X^*$ by the deg-lex order,
i.e., compare two words first by their degrees, then order them
lexicographically when the degrees are equal.
 Clearly, this order is a monomial order. Now we
define the words
$$x_{ji}=x_{j}x_{j-1}...x_{i},$$ where $j>i>1$
and
$x_{j1\varepsilon}=x_{j}x_{j-1}...{x_{1}}^{\varepsilon},\varepsilon=\pm1$.
\\

The proof of the following lemma is straightforward. We omit the
details.
\begin{lemma}\label{l5}
For $\varepsilon=\pm1$, the following relations hold in the
alternating group $A_{n}$:
\end{lemma}
\begin{enumerate}
\item[(4.1)]$x_{1}^{2\varepsilon}=x_{1}^{-\varepsilon}$
\item[(4.2)]$x_{i}^{2}=1,(i>1)$
\item[(4.3)]$x_{j}x_{i}=x_{i}x_{j},(j-1>i\geq2)$
\item[(4.4)]$x_{j}x_{1}^{\varepsilon}=x_{1}^{-\varepsilon}x_{j},(j>2)$
\item[(4.5)]$x_{ji}x_{j}=x_{j-1}x_{ji},(j>i\geq2)$
\item[(4.6)]$x_{j1\varepsilon}x_{j}=x_{j-1}x_{j1-\varepsilon},(j>2)$
\item[(4.7)]$x_{2}{x_{1}}^{\varepsilon}x_{2}=x_{1}^{-\varepsilon}x_{2}x_{1}^{-\varepsilon}$
\item[(4.8)]$x_{1}^{\varepsilon}x_{1}^{-\varepsilon}=1$
\end{enumerate}

Now, we can state the normal form theorem for the group $A_n$.

\begin{theorem} \label{ta2}
Let
\begin{eqnarray*}
A_{n}=gp\langle x_{i} \ (1\leq i\leq n-2)&;&x_{1}^{3}=1, \
(x_{i-1}x_{i})^3={x_{i}}^{2}=1 \ (2\leq i\leq n-2),\\
&&(x_{i}x_{j})^{2}=1 \ (1\leq i<j-1,j\leq n-2)\rangle
\end{eqnarray*}
be the alternating group of degree $n$. Let
$S=\{x_{1}^{2\varepsilon}-x_{1}^{-\varepsilon}, \ x_{i}^{2}-1 \
(i>1), \ x_{j}x_{i}-x_{i}x_{j} \ (j-1>i\geq2), \
x_{j}x_{1}^{\varepsilon}-x_{1}^{-\varepsilon}x_{j} \ (j>2), \
x_{ji}x_{j}-x_{j-1}x_{ji} \ (j>i\geq2), \
x_{j1\varepsilon}x_{j}-x_{j-1}x_{j1-\varepsilon} \ (j>2), \
x_{2}{x_{1}}^{^{\varepsilon}}x_{2}-x_{1}^{-\varepsilon}x_{2}x_{1}^{-\varepsilon},
\ x_{1}^{\varepsilon}x_{1}^{-\varepsilon}-1, \ \varepsilon=\pm1\}$,
where $x_{ji}=x_{j}x_{j-1}...x_{i}, \ j>i>1, \
x_{i1\varepsilon}=x_{i}x_{i-1}...{x_{1}}^{\varepsilon}, \
\varepsilon=\pm1$. Then
\begin{enumerate}
\item[(i)] \ $S$ is a Gr\"{o}bner-Shirshov basis of the alternating group
$A_{n}$;
\item[(ii)] \ every element $w$ of $A_n$ has a unique representation
$w=x_{1j_{1}}x_{2j_{2}}...x_{n-2j_{n-2}}$, where $x_{tt}=x_{t} \
(t>1), \ x_{ii+1}=1, \ x_{i1}=x_i x_{i-1}\cdots x_1^{\varepsilon},
  \ 1\leq j_{i}\leq i+1, \ 1\leq
i\leq n-2, \ \varepsilon=\pm1$ (here we use $x_{i1}$ instead of
$x_{i1\varepsilon}$).
\end{enumerate}
\end{theorem}
\noindent\textbf{Proof } By Lemma \ref{l5}, it is easy to see that
every element $w$ of $A_n$ has a representation
$w=x_{1j_{1}}x_{2j_{2}}...x_{n-2j_{n-2}} \ (j_i\leq i+1)$. Here
$x_{1j_{1}}$ may have 3 possibilities $1,x_1,x^{-1}_1; \ x_{2j_{2}}$
 4 possibilities, and generally, $x_{ij_{i}} \ \ i+2$
possibilities, $1\leq i\leq n-2$. So there are $n!/2$ words. From
this fact, it follows that each representation is unique since
$|A_{n}|=n!/2$. On the other hand, it is clear that $Red(S)$
consists of the same words $w$. Hence, by Composition-Diamond Lemma,
$S$ is a Gr\"{o}bner-Shirshov basis of the alternating group
$A_{n}$. \ \ $\square$
\\

\noindent{\bf Remark}: According to \cite{ls}, normal form in
$S_{n-1}$ is $x_{1j_{1}}x_{2j_{2}}...x_{n-2j_{n-2}}$, but here
$x_{i1}=x_i x_{i-1}\cdots x_1$.
\\

\noindent{\bf Acknowledgement}: The authors would like to express
their deepest gratitude to Professor L. A. Bokut for his kind
guidance, useful discussions and enthusiastic encouragement during
his visit to the South China Normal University.

\ \

\end{document}